\input amstex.tex
\input amsppt.sty   
\magnification 1200
\vsize = 8.5 true in
\hsize=6.2 true in
\NoRunningHeads 
\nologo       
\parskip=\medskipamount
        \lineskip=2pt\baselineskip=18pt\lineskiplimit=0pt
       
        \TagsOnRight
        \NoBlackBoxes

        \topmatter
        \title
        Non-Lipshitz flow of the nonlinear Schr\"odinger equation on surfaces
        \endtitle
\author
         W.-M.~Wang        \endauthor        
\address
{D\'epartement de Math\'ematique, Universit\'e Paris Sud, 91405 Orsay Cedex, FRANCE}
\endaddress
        \email
{wei-min.wang\@math.u-psud.fr}
\endemail
\abstract
We construct non-Lipshitz flow in $H^s$ for  the cubic nonlinear Schr\"odinger equation
on the 2-torus of revolution with a Lipshitz or smooth metric . The non-Lipshitz property holds
for all $s<2/3$ for Lipshitz metric and $s<1/2$ for smooth metric. Both coincide with the Sobolev 
exponents for uniform local well-posedness.  
 \endabstract
\endtopmatter
\document
\head{\bf 1. Introduction}\endhead
We consider the Cauchy problem for the cubic nonlinear Schr\"odinger equation on compact Riemannian surfaces without boundary:
$$
\cases i\frac\partial{\partial t}u =-\Delta u+|u|^{2}u,\\
u(t=0)=u_0, \endcases\tag 1
$$
with either a smooth or Lipshitz metric $g$ and $\Delta$ is the corresponding Laplace-Beltrami operator. 
(Surfaces with smooth boundary are included in Lipshitz $g$.)
We say that the Cauchy problem (1) is uniformly locally well-poesed in $H^s$ if for all $R>0$, there exists $T>0$ and a 
Banach space $X_T$ included in  $\Cal{C}^0((-T,T), H^s)$ such that for all $f \in H^s$ with $\|f\|_{H^s} \leq R$, (1.1) has
a unique solution $u \in X_T$ and the flow map is uniformly continuous in $\Cal C^0$. The problem is globally well-posed
if $T$ can be arbitrary large (without the uniformity requirement).  

On the flat torus, using multiple space-time Fourier series, Bourgain \cite {B} proved that (1) is
locally well-posed in $H^s$ for all $s>0$.   Later  Burq, Gerard and Tzvetkov \cite{BGT1} proved
local well-posedness in $H^s$ for $s>1/2$ for smooth $g$; Blair, Smith and Sogge proved $s>2/3$ for 
Lipshitz $g$ (cf. also \cite{A}). The latter two use dispersive estimates valid on short time intervals.
So (1.1) is energy subcritical and has global solutions in $H^s$ for $s\geq 1$. Recently, Hani \cite{H} 
proved moreover that (1) has global solutions in $H^s$ for $s>2/3$ and smooth $g$ .
 
On the flat 2-torus, it is known \cite{B, CCT} that $s>0$ is optimal. On the 2-sphere, it is known \cite {BGT2, 3} that 
$s>1/4$ is optimal instead of $1/2$. The purpose of this note is to present a simple construction valid for smooth as well 
as Lipshitz $g$ to exhibit non-Lipshitz flow for $s<1/2$, $2/3$ respectively. For Lipshitz $g$, we assume that
the singularities are of type $|x|$.

Concretely,  on the 2-torus of revolution we prove 
\proclaim {Proposition}
 Let $ds^2=dx^2+g(x)dy^2$ with $g\in C^3$ or Lipshitz with isolated singularities and admitting a unique 
 global maximum.  Then there are initial data, which are eigenfunctions of the Laplacian, such that the flow map is not Lipshitz
 in $H^s$ for $s<1/2$ when $g$ is $C^3$ and $s<2/3$ when $g$ is only Lipshitz. 
 \endproclaim
 
\noindent{\it Remark.} We relate the time scale and the Sobolev scale when this non-Lipshitz behavior is manifest.   
The proposition uses data at single high frequency $k\gg 1$ and non-Lipshitz behavior
is observed at time scale $t\asymp k^{1/2}$ for smooth $g$  and $t\asymp k^{2/3}$ for Lipshitz $g$. For $t$ up to $\Cal O(1)$, non-Lipshitz flow is observed 
for $s<1/4$ and $1/3$ respectively. 

We note that this transition from $1/4$ to $1/2$ also occurs on the 2-sphere. So the Sobolev exponent $1/4$ obtained in
\cite{BGT2, 3} is strictly local.  

Previously, for smooth $g$, non-Lipshitz flow is known \cite {T} for $s<1/4$ 
under the assumption of a stable non-degenerate periodic geodesic for $t\to 0$ as a negative power of $k$ using semi-classical
constructions.  The fact that $1/2$ and $2/3$ are observed on the torus of revolution, but at a longer time scale essentially reflects the
stability of high frequency data, cf. \cite{H} in the present context and \cite{W1, 2} (also the review paper \cite{W3}) in the energy supercritical
context.  
\bigskip

\head{\bf 2. Proof}\endhead
 Let the torus be the set $[-\pi, \pi)^2$ with the metric $g$ and identify the end points. From separation of 
 variables, the spectrum of the Laplacian decomposes into:
 $$\sigma(\Delta)=\bigoplus_k\sigma (-\frac {d^2}{dx^2}+k^2 g^{-1}),\tag 2$$ 
 where $k\in\Bbb Z$ is the Fourier dual of $y$. We investigate, in the high frequency limit: $|k|\gg 1$,
 the ground state eigenfunction and the first two eigenvalues of the Schr\"odinger operator defined in the right side of (2). 
 Below we assume $k$ is positive, as negative $k$ works likewise.  
 
 Assume $g$ is only Lipshitz, $g$ smooth works similarly, cf. \cite{C}. Since $g$ has only isolated singularities and
 a unique global maximum and we are interested in the semi-classical limit $k\gg 1$, it suffices to 
 take $g=(|x|+1)^{-1}$ with the unique singularity and global maximum at $x=0$. Other cases can be reduced to this. 
 
 The Schr\"odinger operator is then
 $$H=-\frac {d^2}{dx^2}+k^2 |x|+k^2,\tag 3$$  
 on $\Bbb L^2[-\pi, \pi)$ with periodic boundary conditions. The reference operator is therefore
 $$\Cal A= -\frac {d^2}{dx^2}+k^2 |x|\tag 4$$ on $L^2(\Bbb R)$.
 It is classical that $\Cal A$ has eigenvalues $\alpha_n$, which are deduced from the
 zeroes of the Airy function or its derivative and that its eigenfunctions $\psi_n$ are obtained from
 the Airy function for positive arguments so that the $n$th eigenfunction has parity $(-1)^n$, 
$n=0$, $1$...
 
 More precisely, $\psi_n$ can be written as 
 $$\align \psi_n(x)=&C_n Ai(k^{2/3}|x|-a_n),\qquad \qquad\, x\geq 0,\\
 =&(-1)^nC_n Ai(k^{2/3}|x|-a_n),\qquad x<0,\endalign$$
 where $C_n$ is a nomalizing constant and for $n$ even, $a_n$ is the $n/2+1$ zero of the derivative of 
 the Airy function $Ai'(x)$ and for  $n$ odd, $a_n$ is the $(n+1)/2$ zero of 
 the Airy function $Ai(x)$.  
 
 Below $k$ is large, $k\gg 1$. For us, it suffices to know that 
 $$\alpha_0,\text{ and }\alpha_1-\alpha_0=\Cal O(k^{4/3}),\tag 5$$ 
 $$\Vert \psi_0\Vert_\infty=k^{1/3}\Vert \psi_0\Vert_2, \tag 6$$ 
 and 
 $$|\psi_0(x)|\sim k^{1/6}\frac{e^{-k|x|^{3/2}}}{|x|^{1/4}},\quad |x|>k^{-2/3}.\tag 7$$
 
Let $\phi_0$ be the ground state eigenfunction, $\lambda_0$ and $\lambda_1$ the first two eigenvalues of $H$. 
Using (5, 7) and standard perturbation theory, we then obtain
$$\align &\lambda_0=k^2+\Cal O(k^{4/3})+\Cal O(e^{-k/2})=k^2+\Cal O(k^{4/3})\tag 8\\
&\lambda_1-\lambda_0=\Cal O(k^{4/3}),\tag 9\endalign$$ 
and 
$$\phi_0=\psi_0+\Cal O(e^{-k/2}),\tag 10$$
where the $\Cal O$ is in $L^2$, which in turn gives  
 $$\Vert \phi_0\Vert_\infty  \asymp k^{1/3}\Vert \phi_0\Vert_2\tag 11$$ 
using (6) and Sobolev embedding in one dimension.   
\bigskip
We now proceed to study the Cauchy problem (1) with the initial condition $u_0$:
$$u_0(x, y)=ae^{iky}\phi_0(x),\tag 12$$ where 
$\phi_0$ is assumed to be normalized, $\Vert \phi_0\Vert_2=1$, $a=\Cal O(k^{-s})$, $s\geq 0$, so that  $\Vert u_0\Vert_{H^s}=\Cal O(1)$. 
The solution $u$ can be written as 
$$u(x, y, t)=e^{iky}v(x, t)\tag 13$$
with $v$ satisfying
$$
\cases i\frac\partial{\partial t}v =(-\frac{d^2}{dx^2}+k^2 g^{-1})v+|v|^{2}v,\\
v(t=0)=a\phi_0.\endcases\tag 14
$$

We seek $v$ in the form
$$v(x, t)=a\gamma(t)e^{it(\lambda_0+a^2\omega)}\phi_0(x)+\sum_{j=1}^\infty q_j(t)\phi_j(x),$$
where $\gamma(0)=1$, $q_j(0)=0$, 
$$\omega=\frac {1}{2}\Vert \phi_0\Vert_4^4\tag 15$$ and $\phi_j$ is the $j$th eigenfunction of the 
Schr\"odinger operator in (3). We note that $\omega$ is the frequency modulation that is at the root of
this non-Lipshitz flow.

From energy conservation, we have 
$$a^2\lambda_0+\frac{1}{2}a^4\Vert \phi_0\Vert_4^4=a^2\gamma^2(t)\lambda_0+\sum_{j=1}^\infty \lambda_j|q_j(t)|^2+\frac{1}{2}\Vert v\Vert_4^4,\tag 16$$
since $\Vert \phi_0\Vert_2=1$. Using $L^2$ conservation: 
$$a^2\gamma^2+\sum |q_j|^2:=a^2\gamma^2+\Vert q\Vert_2^2=a^2, \tag 17$$
we then obtain from (16)
$$\aligned \frac{1}{2}a^4\Vert \phi_0\Vert_4^4&=\sum(\lambda_j-\lambda_0) |q_j|^2+\frac{1}{2}\Vert v\Vert_4^4\\
&\geq k^{4/3}\Vert q\Vert_2^2+\frac{1}{2}\Vert v\Vert_4^4,\endaligned\tag 18$$
where we used (9) to reach the last estimate.

(18) gives the following estimates valid for all $t$:
$$\Vert q\Vert_4\leq \Vert v\Vert_4+a\Vert \phi_0\Vert_4\leq 2a\Vert \phi_0\Vert_4,\tag 19$$
$$\Vert q\Vert_2\leq k^{-2/3}a^2\Vert \phi_0\Vert_4^2.\tag 20$$
(17, 20) then give for all $t$: 
$$(1-\gamma^2(t))\leq k^{-4/3}a^2\Vert \phi_0\Vert_4^4.\tag 21$$
Further, the first line of (18) gives 
$$\aligned \frac{1}{2}a^4\Vert \phi_0\Vert_4^4&\geq k^{4/3(1-s)}\sum(\lambda_j-\lambda_0)^s |q_j|^2\\
&\geq k^{4/3(1-s)}(\Vert q\Vert_{H^s}^2-\lambda_0^s\Vert q\Vert_2^2),\quad s<1,\endaligned$$
which leads to
$$\Vert q\Vert_{H^s}\leq a^2\Vert \phi_0\Vert_4^2k^{s-2/3},\quad s<1,\tag 22$$
where we also used (20).

To solve for $\gamma$, we project $v$ onto $\phi_0$ and obtain
$$\gamma(t)=a^{-1}e^{it(\lambda_0+a^2\omega)}(v, \phi_0).$$
So time derivative satisfies
$$\aligned i\dot{\gamma}(t)&=-a^{-1}(\lambda_0+a^2\omega)e^{it(\lambda_0+a^2\omega)}(v, \phi_0)\\
&\qquad + a^{-1}e^{-it(\lambda_0+a^2\omega)}([(-\frac{d^2}{dx^2}+k^2 g^{-1})v+|v|^{2}v], \phi_0),\endaligned$$
where we used (14).
The choice of $\omega$ cancels the leading order nonlinear term and we obtain
$$\aligned |\dot{\gamma}|\leq &a^2\omega (1-|\gamma|^2)|\gamma|\\
&+ a^2\Cal O(\Vert q\Vert_2\Vert\phi_0\Vert_6^3+\Vert\phi_0\Vert_{\infty}^2\Vert q\Vert_2^2|\gamma|+\Vert\phi_0\Vert_{\infty}\Vert q\Vert_4^2
\Vert q\Vert_2).\endaligned$$
Using $$a=\Cal O (k^{-s}),\, \Vert \phi_0\Vert_4\asymp k^{1/6},\,  \Vert \phi_0\Vert_6\asymp k^{2/9}\, \text{and } \Vert \phi_0\Vert_\infty\asymp k^{1/3},$$
we obtain $$\dot{\gamma}=\Cal O (k^{-4s}+k^{-1/3-4s}+k^{-6s}+k^{-6s+1/3}), \, s<1.$$
We note that the precise value of the right side is not important as long as $\dot{\gamma} t$ is sufficiently small.

Let $S_t$ be the flow at $t$ (if it exists). We note that for $u_0$ of the form (12), the solution $u_t$ exists globally in $H^1$. 
Let $s\sim 2/3^-:=2/3-\delta$ for arbitraryly small $\delta>0$. 
Choose two initial data as in (12) with $a_1=k^{-2/3+\delta}$, $a_2=a_1+\epsilon$ with $\epsilon=k^{-2/3-2\delta}$
and $t\asymp k^{2/3}$. 

Let $a=a_1$, $a_2$ and $q=q_1$, $q_2$ be the remainders, then using also (15, 22) we have 
$$\aligned \Vert q\Vert_{H^s}\leq k^{-1+\delta}\ll &a\epsilon\omega t\asymp k^{-\delta}\ll 1,\\
&a^2\omega t\asymp k^{2\delta}\gg 1.\endaligned$$
So
$$\Vert S_t\Vert_{\text{Lip}}\geq \frac{\Vert u_t^{(2)}-u_t^{(1)}\Vert_{H^s}}{{\Vert u_0^{(2)}-u_0^{(1)}\Vert_{H^s}}}\asymp a^2\omega t\asymp k^{2\delta}\gg 1$$
for $t \asymp k^{2/3}$, $k\gg 1$.

Using exactly the same argument for smooth $g$ but Hermite instead of Airy function gives
$$\Vert S_t\Vert_{\text{Lip}}\geq k^{0^{+}}$$
for $t \asymp k^{1/2}$ in $H^{1/2^-}$, cf. \cite{C}.
\hfill $\square$


\bigskip


\Refs\nofrills{References}
\widestnumber\key{CFKSA}

\ref
\key {\bf A}
\by R. Anton
\paper  Strichartz inequalities for Lipshitz metrics on manifolds and the nonlinear Schr\"odinger equation on domains
\jour  Bull. SMF
\vol 136
\pages 27-65
\yr 2008
\endref



\ref
\key {\bf BSS}
\by M. Blair, H. Smith, C. Sogge
\paper  On Strichartz estimates for Schr\"odinger operators in compact manifolds with boundary
\jour Proc. Amer. Math. Soc. 
\vol 138
\pages 247-256
\yr 2008
\endref


\ref
\key {\bf B}
\by J. Bourgain
\paper  Fourier transformation restriction phenomena for certain lattice subsets and applications to
nonlinear evolution equations, part I: Schr\"odinger equations
\jour Geom. and Func. Anal.
\vol 3
\pages 107-156
\yr 1993
\endref


\ref
\key {\bf BGT1}
\by N. Burq, P. Gerard, N. Tzvetkov
\paper  Strichartz inequalities and the nonlinear Schr\"odinger equations on compact manifolds
\jour Amer. J. Math.
\vol 126
\yr  2004
\pages 569-605
\endref

\ref
\key {\bf BGT2}
\by N. Burq, P. Gerard, N. Tzvetkov
\paper  Bilinear eigenfunction estimates and the nonlinear Schr\"odinger equations on surfaces
\jour Invent. Math. 
\vol 159
\yr  2005
\pages 187-223
\endref

\ref
\key {\bf BGT3}
\by N. Burq, P. Gerard, N. Tzvetkov
\paper  An instability property of the nonlinear Schr\"odinger equation on $S^d$
\jour Math. Res. Lett. 
\vol 9
\yr  2002
\pages 323-335
\endref


\ref
\key {\bf C}
\by F. Catoire
\paper \'Equation de Schr\"odinger non-lin\'eaire dans le tore plat g\'en\'erique et le tore de r\'evolution 
\jour Th\'ese Universit\'e Paris-Sud
\yr  2010
\pages 
\endref

\ref
\key {\bf CCT}
\by M. Christ, J. Colliander, T. Tao
\paper  Asymptotics, frequency modulation and low regularity ill-posedness for 
canonical defocusing equations
\jour Amer. J. Math. 
\vol 125
\pages 1235-1293
\yr 2003
\endref


\ref
\key {\bf H}
\by Z. Hani
\paper  Global well-posedness of the cubic nonlinear Schr\"odinger on compact manifolds without boundary
\jour Arxiv: 1008. 2826
\vol 
\pages
\yr 2010
\endref

\ref
\key {\bf T}
\by L. Thomann
\paper  The WKB method and geometric instability for nonlinear Schr\"odinger equations on surfaces
\jour Bull. SMF 
\vol 136
\pages 167-193
\yr 2008
\endref


\ref
\key {\bf W1}
\by W.-M. Wang
\paper  Supercritical nonlinear Schr\"odinger equations I : Quasi-periodic solutions
\jour Arxiv: 1007, 0154
\vol 
\pages 
\yr 2010
\endref

\ref
\key {\bf W2}
\by W.-M. Wang
\paper  Supercritical nonlinear Schr\"odinger equations II : Almost global existence
\jour Arxiv: 1007, 0156
\vol 
\pages 
\yr 2010
\endref

\ref
\key {\bf W3}
\by W.-M. Wang
\paper  Spectral methods in PDE
\jour Milan J. Math., 
\vol 78, no. 2
\pages Arxiv: 1009.0993
\yr 2010
\endref


\endRefs
\enddocument
\end